\documentclass[10pt]{article}
\usepackage{latexsym,amssymb}
\usepackage{enumerate, verbatim}
\usepackage{graphics}

\newtheorem{thm}{\sc Theorem}[section]
\newtheorem{prop}[thm]{\sc Proposition}
\newtheorem{blem}[thm]{\sc Basic Lemma}
\newtheorem{lem}[thm]{\sc Lemma}
\newtheorem{cor}[thm]{\sc Corollary}

\newcommand{\CF}{\medskip{\it Francesca,}\,\,\,}

\newcommand{\G}{{\cal G}}

\newcommand{\Pa}{{\rm P}}

\newcommand{\Sy}{{\rm Sym\,}}
\newcommand{\Al}{{\rm Alt\,}}
\newcommand{\gD}{\Delta}
\newcommand{\ga}{\alpha}
\newcommand{\gb}{\beta}
\newcommand{\gga}{\gamma}
\newcommand{\gG}{\Gamma}

\newcommand{\gO}{\Omega}
\newcommand{\gS}{\Sigma}
\newcommand{\go}{\omega}
\newcommand{\pf}{{\it Proof:\,} }
\newcommand{\dne}{\hfill $\Box$ \bigskip}

\newcommand{\lnr}{{\cal L}_{\rm NR}}
\newcommand{\lur}{{\cal L}_{\rm SR}}
\newcommand{\lst}{{\cal L}_{\rm ST}}
\newcommand{\loe}{{\cal L}_{\rm OE}}

\newcommand{\e}{\!\!:\!}

\begin{document}

\title{\sc  On Orbit Equivalence and Permutation Groups Defined by Unordered Relations }

\author{ Francesca Dalla Volta \\  
{\small Dipartimento di Matematica e Applicazioni, Universit\`{a} Milano Bicocca}\\ 
{ \small 20125 Milano, Italy}\\
{\small\tt francesca.dallavolta@unimib.it}\\\\
Johannes Siemons\\  
{\small School of Mathematics, University of East Anglia,}\\ 
{\small Norwich, NR4 7TJ, United Kingdom}\\{\small\tt j.siemons@uea.ac.uk}}
\date{\scriptsize {Version of 16 October 2010, printed \today}}
\maketitle

\begin{abstract}
\noindent

For a set $\gO$ an {\it unordered relation} on $\gO$ is a family $R$ of
subsets of $\gO.$ If $R$ is such a relation we let $\G(R)$ be the group of all permutations on $\gO$ that preserves $R,$ that is $g$ belongs to $G(R)$ if and only if  $x\in R$  implies $x^{g}\in R.$ We are interested in permutation groups which
can be represented as $G=\G(R)$ for a suitable unordered relation $R$
on $\gO.$  When this is the case, we say that $G$ is defined by the relation $R,$ or that $G$ is a relation group.  We prove that a  primitive permutation group $\neq \Al(\gO)$ and of degree $\geq 11$ is a relation groups. The same is true for many classes of  finite imprimitive groups, and we give general conditions on the size of blocks of imprmitivity, and the groups induced on such blocks, which guarantee that the group is defined by a  relation. 

This property is closely connected  to the orbit closure of permutation groups. Since relation groups are  orbit closed the results here imply  that many classes of imprimitive permutation groups are orbit closed.

\phantom{XXXXXXXXXXX} of set systems

\end{abstract}

\section{Introduction}

Let  $\gO$ is a set. Then the symmetric group $\Sy(\gO)$ acts naturally on the  collection $P(\gO)$ of all subsets of $\gO.$ If $R$ is a subset of $P(\gO)$ we may define the group  
$$\G(R):=\{\,g\in\Sy(\gO)\,\,:\,\,x^{g}\in R \,\,\mbox{ for all}\,\, x\in R\}$$ 
of all permutations that leave $R$ invariant.  In this paper we are interested in those  permutation groups which arise in this way. We consider $R$ as an {\it unordered relation, } or  just a {\it relation} on $\gO,$ and we say that a permutation group $G$ on $\gO$ is  a {\it relation group } if there is a relation $R$ on $\gO$ so that $G=\G(R).$ The notion of relation group has appeared earlier in the literature, see for instance Betten's paper~\cite{betten}. Our first result is Theorem~\ref{thm:PrimRela}\, where we show that a finite primitive group $G$ on $\gO$ is a relation group unless $G$ is $\Al(\gO)$ or one of  ten exceptions of degree $\leq 10.$ The remainder of the  paper therefore deals with finite imprimitive groups. In Corollary~\ref{cor:many}\, it is shown that most imprimitive groups are relation groups, in the following sense:\, If $H$ is an imprimitive group of degree $n$ with a block of imprimitivity $\gD$  so that the group induced by $H$ on $\gD$ does not contain $\Al(\gD)$  and so that $\gD$ is sufficiently large in comparison to $n|\gD|^{-1}$ then all subgroups of $H$ are relation groups.  In Theorem~\ref{inductionGluck} the same conclusion is obtained if some  imprimitivity chain for $H$ does not include factors belonging to a certain class of groups defined there. This allows us also to deal efficiently with solvable groups.

This study of relation groups arose from a closely connected  property of permutation groups.  Two groups on $\gO$ are called  {\it orbit equivalent } if they have the same orbits on $P(\gO).$ From this the {\it orbit closure } of $G$ is defined as the largest group $G^{*}$ on $\gO$ that is orbit equivalent to $G,$ and $G$ is {\it orbit closed } if $G=G^{*}.$ Orbit equivalence  has been studied for some time, and we give some  references in the text.  From the classification of finite simple groups there is a complete list of primitive pairs of orbit equivalent groups due to Seress~\cite{Seress1} based on an observation by Cameron, Neumann and Saxl~\cite{CNS}. More recently  Seress and Young~\cite{Seress2}\, investigate the case when one of the groups is the wreath product of two primitive groups.  

It is almost immediate that a relation group is orbit closed, see Proposition~\ref{2.5}. The results mentioned above therefore show that large classes of permutation groups are orbit closed.   In fact, we show that a primitive group  is orbit closed if and only if it is a relation group. This however fails for infinitely many imprimitive groups, see Corollary~\ref{cor:WreathProduct1}.

Independently of this application relation groups are interesting in their own right: They include the full automorphism groups of graphs, 
designs, geometries and undirected set systems generally.  Our main tool for dealing with relation groups are regular sets for primitive groups, and this requires the classification of finite simple groups. We also require the explicit knowledge of primitive groups without regular sets in~\cite{Seress1}. In the last section we mention open problems for relation groups and orbit closure.

\bigskip
\section{Prerequisites}

Throughout $\gO$ denotes a set of size $n$ and $P(\gO)$ denotes the collection of all subsets of $\gO.$ For $0\leq k\leq n$  the $k$-element subsets of $\gO$ are denoted  by $\gO^{\{k\}}.$  All sets and groups in this paper are finite. 

\bigskip
{\bf 1. Action on Subsets:}\quad The symmetric group $\Sy(\gO)$ acts naturally   on $\Pa(\gO)$ via $g\!:x\mapsto x^{g}:=\{\go^{g}\,:\,\go\in x\}$ for  $g\in \Sy(\gO)$ and $x\in \Pa(\gO).$ 
Let $(G,\gO)$ be a permutation group on $\gO$  and let $x\in \Pa(\gO).$  Then  $x^{G}=\{\,x^{g}\,:\,g\in G\,\}$ is the  $G$-orbit of $x.$ 
We say that $G$ is {\it set-transitive } if $G$ acts transitively on $\gO^{\{k\}}$ for all $k=0,\!..,n.$
We use standard {\sc Atlas} notation for groups, in particular $C_{m},$ $D_{m}$ and $F_{m}$ are the cyclic, dihedral and Frobenius groups of order $m,$ respectively. We write also $S_{m}=\Sy(\gG)$ and  $A_{m}=\Al(\gG)$ when $|\gG|=m.$

\begin{lem}\label{2.1}
{\rm (Beaumont and Peterson~\cite{Beau}, 1955)}\quad Let $H\subseteq \Sy(\gO)$ be a set-transitive  group  of degree $n\geq 2.$ Then $H\supseteq \Al(\gO)$ or \vspace{-5mm} \begin{enumerate}[(i)]
\item $n=5$ and $H=F_{20},$\vspace{-3mm}
\item $n=6$ and $H=L_{2}(5).2\simeq S_{5}$ or  \vspace{-3mm}
\item $n=9$ and $H\in\{L_{2}(8)\subseteq L_{2}(8).3\}.$ \vspace{-3mm}
\end{enumerate}
\end{lem}

\medskip
{\bf 2. Regular Sets:}\quad Let $(G,\gO)$ be a permutation group on $\gO.$ We denote  the setwise stabilizer of $x$ in $G$ by $G_{x}.$ Then $x$ is a {\it regular set} for $G$ if $G_{x}=1.$ A key tool for this paper is the  characterization  of  primitive groups without a regular set. This is obtained from  the classification of finite simple groups.   We state this result in full as we need the details later on. The set of primitive groups of degree $5\leq n\leq 32$ in the theorem without regular set will be denoted by $\lnr.$

\begin{thm} \label{Seress1}{\rm (Seress~\cite{Seress1}
)} 
\quad Let $H\subseteq \Sy(\gO)$ be a primitive group  of degree $n\geq 2.$ Suppose that $H$ is not set-transitive.  Then $H$ has a regular set  on $\gO$ if and only if $H$ is not one of the following:
\vspace{-5mm} \begin{enumerate}[(i)]
\item $n=5$ with $H=D_{10}\,;$\vspace{-3mm}
\item $n=6$ with $H=L_{2}(5)\,;$\vspace{-3mm}
\item $n=7$ with  $H=F_{42},$\, $H=L_{3}(2)\,;$   \vspace{-3mm}
\item $n=8$ with 
             $H\in \{2^{3}\e 7.3\subset 2^{3}\e L_{3}(2)\},$\,$H\in\{L_{2}(7)\subset
L_{2}(7).2\}\,;$\vspace{-3mm}
\item $n=9$ with  $H\in\{3^{2}\e (2.L_{2}(3))\subset 3^{2}\e 2.L_{2}(3).2\},$\,
$H\in\{3^{2}\e D_{8}\subset 3^{2}\e 8.2\subset 3^{2}\e 2.L_{2}(3).2\}\,;$ \vspace{-3mm}
\item $n=10$ with  $H\in\{S_{5},\,L_{2}(9)\subset L_{2}(9).2\subset L_{2}(9).2.2\}\,$ \\$H\in\{L_{2}(9)\simeq PSU(2,9)\subset PGU(2,9)\},\,$ $H\in\{L_{2}(9)\subset M_{10}\}\,;$ \,\vspace{-3mm}
\item $n=11$ with  $H\in\{L_{2}(11)\subset M_{11}\}\,;$  \vspace{-3mm}
\item $n=12$ with $H=L_{2}(11).2,\,$ $H\in\{M_{11}\subset M_{12}\}   \,;$   \vspace{-3mm}
\item $n=13$ with  $H=L_{3}(3)\,;$   \vspace{-3mm}
\item $n=14$ with  $H=L_{2}(13).2\,;$   \vspace{-3mm}
\item $n=15$ with  $H=L_{4}(2)\simeq A_{8} \,;$   \vspace{-3mm}
\item $n=16$ with  $H=2^{4}\e (A_{5}\times 3).2,$\, 
$H\in \{2^{4}\e A_{6}\subset 2^{4}\e S_{6}\subset 2^{4}\e L_{4}(2)\},$\, \\
$H\in \{2^{4}\e A_{7}\subset 2^{4}\e L_{4}(2)\}\,;$   \vspace{-3mm}
\item $n=17$ with  $H\in \{L_{2}(16).2\subseteq L_{2}(16).4\}\,;$   \vspace{-3mm}
\item $n=21$ with  $H=L_{3}(4).3.2\,;$   \vspace{-3mm}
\item $n=22$ with  $H\in \{M_{22}\subseteq M_{22}.2\}\,;$   \vspace{-3mm}
\item $n=23$ with  $H=M_{23}\,;$   \vspace{-3mm}
\item $n=24$ with  $H=M_{24}\,;$   \vspace{-3mm}
\item $n=32$ with  $H=2^{5}\e L_{5}(2)\,.$   \vspace{-2mm}
\end{enumerate} \vspace{-3mm}
In each case the  containments of the primitive  exceptions are as indicated.  In particular, a primitive, not set-transitive group of degree  $n\not\in\{5,\,6,\,7,\,8,\,9,\,10,\,11,\,12,\,13,\,14,\,15,\,16,\,17,\,21,\,22,\,23,\,24,\,32\}$  has a regular set.
 \end{thm}

{\sc Comments:}\, 1. The theorem requires the classification of finite simple groups in the following way: If $(G,\gO)$ is a permutation group without regular sets then each  of the $2^{|\gO|}$ subsets of $\gO$ is stabilized by a non-identity element. Since a permutation composed of $c$ cycles stabilizes exactly $2^{c}$ subsets, a bound for the minimum degree in $G$ (the least number of elements by a non-identity element in $G$) implies that $G$ is  `large'. On the other hand, from the O'Nan-Scott Theorem and the classification  of finite simple groups an upper bound for the order of  a primitive group not containing the alternating group is obtained, and this yields an upper bound for the order of $G.$ This argument due to Cameron, Neumann and Saxl~\cite{CNS} shows that all but finitely many primitive groups not containing the alternating group of the same degree have a regular set.  In Section 7 we use the same idea to prove that all but finitely many such  groups have regular sets of at least two different cardinalities. The  list $\lnr$ of exceptions is obtained   in Seress' theorem by further careful analysis of the bounds and by direct  computations. 

2.  There are direct methods to establish the existence of a regular set for many classes of permutation groups, and these methods do not  require primitivity or the classification. If the group is given as the automorphism group of a geometrical or combinatorial  object, say a graph, an affine or projective space, or a design more generally, it is often possible to construct a regular set from  the combinatorics of the structure.  In the literature on graphs and set systems  the more general notion  of a {\it distinguishing partition} appears, introduced by Albertson and Collins~\cite{Alber}. These are partitions into an arbitrary number of classes so that only the identity  automorphism stabilizes all  classes, see the paper~\cite{Lafla} by C. Laflamme, Nguyen Van Th\'{e} and N. Sauer. Regular sets therefore correspond to distinguishing partitions into two classes.

3. 
We mention some papers where regular sets are obtained by such direct  methods.  Gluck~\cite{gluck} determines all solvable primitive groups without a regular set. In   \cite{key3, key4, volta1, volta2}  it is shown that regular sets exist for the  affine, projective, unitary or  orthogonal groups in their natural action, apart from a small number of explicitly listed exceptions. Often the size of a smallest regular set is determined as well. 
The  problem of determining the primitive groups in which for given $k$ {\it all} $k$-element subsets are   regular is considered in Bates, Bundy, Hart and Rowley~\cite{bates}.

\bigskip
{\bf 3. Relations and Relation Groups:}\quad We consider  a subset $R$ of $P(\gO)$ as an  {\it unordered relation } on $\gO,$ or just a {\it relation}. (We we shall not consider  the usual ordered relations.)   If $R$ is a relation then 
$$\G(R):=\{\,g\in\Sy(\gO)\,\,:\,\,x^{g}\in R \,\,\mbox{ for all}\,\, x\in R\}$$ is the {\it invariance group\,} of $R,$ or  the  {\it group defined} by $R.$ 
Evidently a permutation group $G$ on $\gO$ is contained in $\G(R)$ if and only if  $R$ is a   union of $G$-orbits on $P(\gO).$ We call  $G$ a {\it relation group}  on $\gO$ if there is some $R$ on $\gO$ for which $G=\G(R). $ It is also possible  to  view $R$ as a {\it hypergraph} with vertex set $\gO,$ and  $\G(R)$ as its full automorphism group. In this paper we prefer however the language of relations as it is natural and more flexible.

Key examples of relation groups arise from combinatorial structures and geometry. For instance, if $\gO$ is the set of  points of a projective or affine space ${\cal S}$ and if $C\subset \gO^{\{3\}}$ is the collinearity relation of  ${\cal S}$ (i.e. the collinear triples) then ${\rm Aut}({\cal S})=\G(C)$ is a relation group by definition. In fact, we note that most classical groups arise in this fashion. 
With this in mind Betten~\cite{betten} calls a relation group also a {\it geometric group.} 
 The problem  of representing a permutation group by an unordered relation on the same set  was also considered in~\cite{key4}. 

\medskip
If $R$ is a relation then its {\it arity}  is the set ${\rm ar}(R)=\{|x|\,\,:\,\,x\in R\}.$  
The arity of  a relation is important for our purpose. The following obvious fact will be used without further mention. 

\begin{lem}\label{2.2}
\quad If $R'$ and $R''$ are relations with ${\rm ar}(R')\cap {\rm ar}(R'')=\emptyset$ then $\G(R'\cup R'')=\G(R')\cap \G(R'').$
\end{lem}

We call  $R$  {\it trivial} if  $R=\bigcup_{k\in I}\,\gO^{\{k\}}$ for some index set  $I\subseteq \{0,\,1,\,..,n\}.$ Evidently $G$ is set-transitive if and only if the only relations $R$ with $G\subseteq \G(R)$ are the trivial relations.

\bigskip
\bigskip
{\bf 4. Orbit Equivalence:}\quad Let $0\leq k\leq n.$ Two permutation groups $G$ and $H$ on $\gO$ are  {\it $k$-orbit equivalent }  if they  have the same orbits on $\gO^{\{k\}}.$ We denote this by $G\approx^{k} H.$ The largest permutation group on $\gO$ that is  $k$-orbit equivalent to $G$ is the {\it $k$-closure} of $G,$ denoted $G^{\{k\}},$ and  $G$ is  {\it $k$-orbit closed, } or just {\it $k$-closed } if $G=G^{\{k\}}.$ Evidently $G^{\{k\}}=G^{\{n-k\}}.$ 
The following theorem states the relationship between the closure groups:

\medskip
\begin{thm} \label{Sie}
{\rm (Siemons~\cite{siem1})}\quad Let $G$ and $H$ be  permutation groups on $\gO$  and suppose that $0\leq k\leq \ell$ satisfy $k+\ell\leq |\gO|.$ Then $G\approx^{\ell} H$ implies $G\approx^{k} H.$ In particular, $G^{\{\ell\}}\subseteq G^{\{k\}}$ and so a  $k$-closed group is $\ell$-closed.
\end{thm}

We say that $G$ and $H$ are  {\it orbit equivalent,\,} denoted $G\approx H,$ if $G\approx^{k} H$ for all $0\leq k\leq n.$ The  {\it orbit closure } of $G$ is the largest permutation group $G^{*}$ on $\gO$ that is  orbit equivalent to $G.$ We say that $G$ is {\it orbit closed } if $G=G^{*}.$ From the theorem it follows that $G^{*}=G^{\{n^{*}\}}$ where $n^{*}=\lfloor \frac n2\rfloor.$ 
There are other closure operations for permutations groups in the literature. These  are derived from  the group action on the cartesian product $\gO^{k},$ rather than $\gO^{\{k\}},$ and were first explored  in  Wielandt's Ohio lecture notes~\cite{Ohio}. While  there are interesting common features,  we emphasize that in this paper we are not concerned with this notion of closure.

The connection between orbit closure and relation groups, of key importance in this paper, is immediate. We have


\begin{prop} \label{2.5}
\quad Let $G$ be a permutation group on $\gO$  and let $R$ be a relation with $G\subseteq \G(R).$ Suppose that $1\leq k\leq \frac {n+1}2$ is an integer so that no set $x$ in $R$ has cardinality $k<|x|<n-k.$ Then $G^{\{k\}}\subseteq \G(R).$ In particular, a relation group is orbit closed. \end{prop}

\pf By Theorem~\ref{Sie}\, the groups $G$ and $G^{\{k\}}$ have the same orbits on $t$-element subsets for all $t\leq k$ and $n-k\leq t.$
Since $G\subseteq\G(R)$ if and only if $R$ is a union of $G$-orbits on $P(\gO)$ we have  
 $G^{\{k\}}\subseteq\G(R)$ provided that no set $x$ in $R$ has cardinality  $k<|x|<n-k.$ \dne

{\sc Comments:}\, 1. For small $k$ there are in some sense `only few' partitions of $\gO^{\{k\}}$ that arise as the orbits of a group on $\gO,$ and therefore `only few' groups that are $k$-orbit closed. To be $k$-closed for small $k$ is therefore a strong property of permutation groups. For instance, $G$ is $1$-closed if and only if $G=\Sy(\gO_{1})\times ...\times \Sy(\gO_{t})$ with   $\gO=\gO_{1}\cup...\cup\gO_{t},$ and $G$ is $2$-closed if and only if there is a (vertex and edge) coloured undirected graph on $\gO$ so that $G$ is its full automorphism group. Much the same can be said  about relation groups $\G(R)$ where the arity of $R$ is bounded by $k$ as in the proposition. So also in this case, to be the group of  such a bounded relation is a strong property of permutation groups. 

 2. There are  transitive orbit closed groups that are not relation groups. The smallest is  the Klein $4$-group $V$ on $\gO=\{1,2,3,4\}.$ It is transitive on sets of size $\neq 2$ and has $3$ orbits on $2$-sets. It is easy to check that $V$  is  orbit closed. Each of the three orbits on $2$-sets gives a relation, and together with their complements these are the only non-trivial relations preserved by $V.$ However, for any choice of an orbit (or its complement) there is a  transposition not in $V$ that  stabilizes the orbit. Hence $V\subset \G(R)$ for any relation preserved by $V.$   In Section 5 we shall show that there are infinitely  many imprimitive groups of this kind. 

3. In Section 4 (Corollary~\ref{cor:PrimRela}) we shall prove that a primitive group is a relation group if and only if it is orbit closed. It is an open problem to formulate general conditions for an orbit closed imprimitive group to be a relation group, see also Corollary~\ref{cor:WreathProduct11}.

\medskip

The primitive orbit closed groups have been determined 
from Lemma~\ref{2.1}\, and 
Theorem~\ref{Seress1}:

\begin{thm} \label{orbeq}
{\rm (Seress \cite{Seress1}
)}\, Let $G\subset H$ be primitive and orbit equivalent permutation groups on $\gO.$ Assume that   $\{G,\,H\}\neq\{\Al(\gO),\,\Sy(\gO)\}.$ Then  \vspace{-5mm} \begin{enumerate}[(i)]
\item $|\gO|=5$ when $\{G,\,H\}=\{C_{5},\,D_{10}\}$ or $\{G,\,H\}\subseteq \big\{F_{20},\,A_{5},\, S_{5}\big\},$ \vspace{-3mm}
\item $|\gO|=6$ when $\{G,\,H\}\subseteq \big\{L_{2}(5).2\simeq S_{5},\,A_{6},\,S_{6}\big\},$\vspace{-3mm}
\item $|\gO|=8$ when $\{G,\,H\}\subseteq \big\{2^{3}\e 7\,\subset\, 2^{3}\e 7.3\,\subset\, 2^{3}\e L_{3}(2)\big\},$\vspace{-3mm}
\item $|\gO|=9$ when $\{G,\,H\}= \big\{3^{2}\e 8\,\subset \,3^{2}\e 8.2\big\},$\, $\{G,\,H\}=\big\{3^{2}\e 2.L_{2}(3)\,\subset \,3^{2}\e 2.L_{2}(3).2 \big\}$ or\\ $\{G,\,H\}\subseteq\big\{L_{2}(8),\, L_{2}(8).3,\,A_{9},\,S_{9}\big\},$\vspace{-3mm}
\item $|\gO|=10$ when $\{G,\,H\}=\big\{L_{2}(9).2\,\subset\, L_{2}(9).2.2\big\}.$   \end{enumerate}
\end{thm}

\bigskip

\section{Subgroups of Relation groups}

Under suitable conditions all subgroups of a relation group are again relation groups. This important and easy to prove property is the subject of this section.  Recall that if $R$ is a relation on $\gO$ then ${\rm ar}(R)=\{|x|\,\,:\,\,x\in R\}$ is the arity of $R.$ 

\medskip
\begin{blem} \label{blem1} 
\, Let $H$ be a permutation group on $\gO$ and suppose that $H$ has a regular set $w.$ Suppose that one of the following holds: \vspace{-4mm} \begin{enumerate}[(i)]
\item There is a relation $R$ on $\gO$ with $H=\G(R)$ so that  $|w|$ does not belong to ${\rm ar}(R),$ or \vspace{-2.5mm}
\item $H$ is not set-transitive, and is a maximal subgroup of $\Sy(\gO)$ with respect to this property. \vspace{-4.5mm}\end{enumerate}
Let $G$ be a subgroup of $H.$ Then $G$  is a relation group on $\gO$ and in particular, $G$ is orbit closed.  
\end{blem}

In this paper we are interested in showing that certain classes of permutation groups are relation groups, and hence are orbit closed. This is an important consequence of being a relation group. However, it makes no sense to repeat this on each occasion. In addition, with Proposition~\ref{2.5} in mind,  even stronger closure properties hold for the group  as  soon as additional information about the arity of the  relation in question in the Basic Lemma is available. 

\medskip\pf (i) For  $G\subseteq H$ let $R^{G}=w^{G}$ be the $G$-orbit of $w$ and let $r$ be the size of $w.$ Then $G\subseteq \G(R\cup R^{G})$ and since $r\not \in {\rm ar}(R)$ we have $\G(R\cup R^{G})=\G(R)\cap \G(R^{G})= H\cap \G(R^{G})$ by Lemma~\ref{2.2}\,. Let  $h$ be an element of $H\cap \G(R^{G}).$ Then there is some $g\in G$ with $w^{h}=w^{g}.$ Hence $hg^{-1}\in H$ stabilizes $w$ so that  $hg^{-1}=1$ and so $h=g$ belongs to $G.$ Therefore $G=\G(R\cup R^{G}).$

(ii) As $H$ is not set-transitive let $x$ be a set with $x^{H}\neq \gO^{\{k\}}$ where $k:=|x|.$ 
Then $H\subseteq \G(x^{H})$ and the latter is not set-transitive. By the maximality condition, $H=\G(x^{H}).$ If $x$ can be chosen so that $k\neq r:=|w|$ the result follows from (i). It remains to consider the case when $H$ is transitive on $\ell$-sets for all $\ell\neq r$ but intransitive on $r$-sets. For any permutation group on $\gO,$ by the Livingstone-Wagner theorem, or by Theorem~\ref{Sie}, transitivity on $t$-sets implies transitivity on $s$-sets provided $s+t\leq |\gO|.$ So in our case $\gO$ must have size $2r.$ As $H$ is transitive on $(r-1)$-sets we have $|H|={2r\choose r-1}\cdot|H_{v}|$ for any set $v$ of size $r-1.$ If $|H_{v}|=1$ we may replace $w$ by $v,$ and the result follows again from (i). Hence  we may assume that $2{2r\choose r-1}\leq |H|.$ But we also have $|H|=|w^{H}|<{2r\choose r}$ since $H$ is intransitive on $r$-element subsets. This implies $r=0,$ a contradiction. 
Finally, if $G$ is a relation groups then $G$ is orbit closed by Proposition~\ref{2.5}.\dne

\medskip
{\sc Comments:}\, 1. This lemma is useful in many situations, we give some examples.  Let $H$ be the group of all collineations of an affine or projective space of dimension $d$ over the field $F,$ viewed as a permutation group on the set $\gO$ of the points of the space. Let $C\subset \gO^{\{3\}}$ be the collinearity relation. Then $H=\G(C),\!$  by definition. As mentioned  in the comments  following Theorem~\ref{Seress1},\, a regular set $w$ for $H$ can be constructed by elementary means, and typically $w$ has cardinality $r=2d$ or $r=2d\pm 1.$ There are some exceptions for small $d$ and $|F|,$ the details are given in   Theorems 3.1 and 4.1 of \cite{key3}. From this information we conclude that any subgroup $G$ of $H$ can be expressed as $G=\G(C\cup w^{G}).$ This relation  is quite simple: only two set sizes are needed. It is evident that $C$ can be replaced by  other  homogeneous relations $C^{*}$ which encapsulate collinearity. We note that it is even possible to find a relation  $C^{*}\subset \gO^{\{r\}}$ so that $\G(C^{*}\cup w^{G})=\G(C^{*})\cap \G(w^{G}),$ see \cite{key4}. For this choice $G=\G(C^{*}\cup w^{G})$ is the group of  a homogeneous relation.  

This example shows that all subgroups of $A\Gamma L(d, F),$ and all subgroups of  $P\Gamma L(d, F),$ are $r$-closed in their action on the points of the  corresponding affine or projective space, apart from a few exceptions for small dimension $d$ and small field $F.$ This follows from Proposition~\ref{2.5}, and by Theorem~\ref{Sie} we may take $r=2d+1.$ The small exceptions are covered in Theorem~\ref{thm:PrimRela}.

\medskip\,
2. It is essential to control the arity of the relations that define a given group. Especially important is the case $|{\rm ar}(R)|=1,$  when $R$ is  homogenous.
In a slightly different setting groups defined by  a homogeneous relation appear as  {\it $2$-representable groups } in Kisielewicz's paper~\cite{Kisie}.
In the second part of the Basic Lemma $H$ is the group of a homogeneous relation, for if $v^{H}$ is any non-trivial $H$-orbit then $H=\G(v^{H}).$ If $G\subset H$ then  the proof has shown that $G=\G(R)$ for a relation  $R$ with $|{\rm ar}(R)|= 2.$

\bigskip
\section{Primitive Groups}

In this section we consider primitive groups. As a  first application of the Basic Lemma we prove that they all are relation groups apart from a few exceptions.

\begin{thm}\label{thm:PrimRelaZero} 
\,Let $H$ be a primitive permutation group on $\gO$ not containing $\Al(\gO).$ Suppose that $|\gO|>9$ and that $H$ does not belong to the list $\lnr$ of  Theorem~\ref{Seress1}.  Then every subgroup of $H$ is a relation group. \end{thm}

{\sc Remark:\,} In the remaining case,  when $H$ does belong to $\lnr,$  it is straightforward to determine all relation groups contained in $H$ by direct computation. However, we  omit these details  here.

\pf As $|\gO|>9$ it follows from Lemma~\ref{2.1}\, that $H$ is not set-transitive and as $H$ is not in $\lnr$  it follows   that $H$ has a regular set. Without loss of generality we may suppose that $H$ is  maximal with the property of not being set-transitive. Now apply the second part of  the Basic Lemma. \dne

Our next aim is to determine the primitive permutation groups on $\gO$ that are not the group of a relation on $\gO.$  We describe these exceptions.  In our lists $(n,H)$ means that $H$ is a group of degree $n=|\gO|.$ 

\medskip
Set-transitive groups: \, If $H\neq \Sy(\gO)$ is  set-transitive group then $H$  is not a relation group.  Apart from $H=\Al(\gO)$ by Lemma~\ref{2.1} the only groups of this type are  
 \vspace{-3mm} 
\begin{enumerate}[\qquad]
\item $\lst=\big\{(5,F_{20}),\,\,(6,\,L_{2}(5).2),\,\,(9,\,L_{2}(8)),\,\,(9,L_{2}(8).3)\big\}\,.$ 
\end{enumerate}
\vspace{-3mm}
 Note that none of these groups  have a  regular set.

\medskip
Other orbit equivalent primitive groups: \,Suppose that $H_{1}\subset H_{2}\neq H_{1}$ are orbit equivalent, and not set-transitive. Then $H_{1}$ is not a relation group, by Proposition~\ref{2.5}. From Theorem~\ref{orbeq}  we obtain the following options for $H_{1}:$
\vspace{-3mm} \begin{enumerate}[\qquad]
\item $\loe=\big\{     (5,C_{5}),\,\,  (8,2^{3}\e 7), \,\,(8,2^{3}\e 7.3),\,\,
(9, 3^{2}\e 8),\,\,(9,3^{2}\e 2.L_{2}(3)),\,\,(10,L_{2}(9).2)   \big\}.$   
\end{enumerate}

\vspace{-3mm}
The groups  $(8,2^{3}\e 7.3),$ $(9,3^{2}\e 2.L_{2}(3))$ and $(10,L_{2}(9).2)$ have no regular sets while the remaining three groups have regular sets.

\medskip
Let ${\cal L}=\,\lst\, \cup\, \loe.$ None of these groups are relation groups, and  together with the alternating groups this determines all groups which  are not relation groups:

\begin{thm}\label{thm:PrimRela} 
 \,Let $H$ be a primitive permutation group on $\gO$ not containing $\Al(\gO).$ Then   $H$ is a relation group if and only if $H$ does not belong to  ${\cal L}.$   \end{thm}

\begin{cor}\label{cor:PrimRela} 
\,A primitive permutation group $H$ on $\gO$ is orbit closed if and only if  it is a relation group. If $H$ is not a relation group then  $H=\Al(\gO)$ or $|\gO|\leq 10.$ In all these cases the  index of $H$ in its orbit closure is divisible by $2.$
\end{cor}

{\sc Remark:}\, Note that primitivity is required. As mentioned before, the Klein $4$-group is orbit closed on $\{1,2,3,4\}$ but not a relation group. We shall discuss other examples in the next sections.

\medskip
{\it Proof of the Corollary:} \,If $H$ is a relation group then  $H$ is orbit closed by Proposition~\ref{2.5}. Conversely, by the comment above, if $H$ is orbit closed then it is not in ${\cal L}$ and so by Theorem~\ref{thm:PrimRela}\, we have that  $H$ is not a relation group. If $H$ is not orbit closed then by the theorem, either $H=\Al(\gO)$ or $H$ belongs to ${\cal L}$. In the latter case it is easy to check that the index of $G$ in its orbit closure is divisible by $2,$ for this see again Theorem~\ref{orbeq}. \dne

\medskip
{\it Proof of Theorem~\ref{thm:PrimRela}:} \,We have already shown that if $H$ belongs to ${\cal L}$ then $H$ is not orbit closed and hence not a relation group. Conversely  suppose that $H$ is primitive on $\gO,$ does not contain $\Al(\gO)$ and does not appear in ${\cal L}.$ Then $H$ is not set-transitive.   Let $X$ be a group containing $H$ that is maximal subject to being not set-transitive.   By maximality $X$ is a relation group. If $X$ has a regular set then the result follows from the second part of the Basic Lemma. 

It remains to inspect the cases when $X$ has no regular set. As $X$  is primitive it is one of the maximal groups in the list $\lnr$ of exceptions in  Theorem~\ref{Seress1}.  (We have $H\subset X\neq H$ but  $H$ may or may not have a regular set.) Using  the {\sc GAP} libraries of primitive groups one can now inspect the primitive subgroups of $X$ which do not belong to ${\cal L}.$ This yields a list of about 58  subgroups $H.$ For all of these one constructs a relation `by hand'.
 In some cases there is a relation group $Y\supset H$ which does have a regular set, and in this case the Basic Lemma suffices to show that also $H$ is a relation group. In the remaining cases we have computed sufficiently many $H$-orbits on $P(\gO)$ so that a relation $R$ can be stated explicitly for which $H=\G(R).$ The following records  the rough details of the computations.  
In all  cases the primitive subgroup $H$ of the group $X\in \lnr$   is a relation group or one of the exceptions in the list ${\cal L}$ of the  theorem.

\smallskip 
{\small
\begin{itemize}
\item Degree $n= 5$ when $H=C_5\subset X=D_{10}$ and $H\in \loe.$
\item Degree $n=6$ when $ X=L_2(5)$ has no primitive proper subgroups.
\item Degree $n=7.$ Here $X_{1}=F_{42}=AGL(1,7)$ (maximal in $S_7)$ or $X_{2}=L_3(2)$ (maximal in $A_7.)$ Both groups are $2$-transitive but neither are $3$- or $4$-homogeneous. Therefore we may represent $X_{1}=\G(R_{1})$ with $R_{1}\subset \gO^{\{3\}}$ and $X_{2}=\G(R_{2})$ with $R_{2}\subset \gO^{\{4\}}.$ One primitive $2$-homogeneous  group appears as  $H_{1}=X_{1}\cap X_{2}=7\e 3=\G(R_{1}\cup R_{2}).$ Apart from this there are only  two other groups, $H_{2}=7\subset H_{1}$ and $H_{3}=7\e 2\subset X_{1}.$ They are $2$-orbit equivalent, with  three orbits of  length $7$ on $\gO^{\{2\}}.$ Any one of them provides a relations $R\subset \gO^{\{2\}}$ with $H_{3}=\G(R).$ Finally, $H_{2}=H_{1}\cap H_{3}=\G(R)\cap \G(R_{1}\cup R_{2})=\G(R\cup R_{1}\cup R_{2}).$

\item Degree  $n=8.$ Here $H_{1}=L_{2}(7)\subset X_{1}=L_{2}(7).2$ or various subgroups of $X_{2}=2^{3}\e L_3(2).$  
Both  groups are $3$-transitive but not $4$-homogeneous. They may be  represented as $X_{1}=\G(R_{1})$  and $X_{2}=\G(R_{2})$  for some $R_{i}\subset \gO^{\{4\}}.$ For $H_{1}$ we take a set $w$ so that $R:=w^{H_{1}}\neq w^{X_{1}},$ see Theorem~\ref{orbeq}. Then $H_{1}=\G(R).$
There are two  primitive subgroups of $X_{2},$  both listed  in $\loe.$ 

\item Degree  $n=9.$ \,There are the two set-transitive groups $L_{2}(8)\subset X_{1}=L_{2}(8).3$  which both appear in  $\lst.$ Apart from this we have $X_{2}=3^{2}\e 2.L_{2}(3).2=AGL(2,3).$ This group has the  primitive subgroups $H_{1}=3^{2}\e4,$ $H_{2}=3^{2}\e Q_{8},$ $H_{3}=ASL(2,3)=3^{2}\e2.L_{2}(3),$ $H_{4}=3^{2}\e D_{8},$ $H_{5}=A\Gamma L(1,9)=3^{2}\e8.3$ and $H_{6}=AGL(1,9).$ These fall into the chains $H_{1}\subset H_{2}\subset H_{3}\subset X_{2},$ further $H_{1}\subset H_{4}\subset H_{5}$ and  $H_{1}\subset H_{6}\subset H_{5}\subset X_{2}.$ The groups $H_{3}$ and $H_{6}$ appear in $\loe,$ they are orbit equivalent to $H_{5}$ and $X_{2},$ respectively.  

It remains to consider $H_{1}\subset H_{2}$ and $H_{4} \subset  H_{5}.$ From computation    $H_{4} \subset  H_{5}$ are not orbit equivalent and so we may take  a set $w_{4}$ with  $R_{4}:=w_{4}^{H_{4}}\neq w_{4}^{H_{5}}.$ Similarly,    $H_{5} \subset  X_{2}$ are not orbit equivalent which allows us to take a set $w_{5}$ with $R_{5}:=w_{5}^{H_{5}}\neq w_{5}^{X_{2}}.$ It is possible to take $|w_{4}|=3$ and $|w_{5}|=4$ with $|R_{5}|=36.$  Then $\G(R_{5})\supseteq H_{5}$ and  $\G(R_{5})\not\supseteq X_{2}$ implies that $\G(R_{5})=H_{5}$ as $X_{2}$ has no orbits of length $36$ on $\gO^{\{4\}}.$  
For a similar reason $\G(R_{4}\cup R_{5}) =H_{4}.$ Next, as $H_{2}$ is not $2$-homogeneous there is $R_{2}\subset \gO^{\{2\}}$ so that  $H_{2}=\G(R_{2}).$ Finally  $H_{1}=H_{4}\cap H_{2}=\G(R_{4}\cup R_{5}\cup R_{2}).$

\item Degree  $n=10.$ \,Here only the possibility $X=P\Gamma L(2,9)$ appears, with $PGL(2,9)$ orbit equivalent to it. Let $H_{1}={\rm Alt}_{5}$ and  $H_{2}={\rm Sym}_{5}$ acting on $2$-sets of a $5$-set, $H_{3}=PSL(2,9),$ $H_{4}=M(10)$ and   $H_{5}=P\Sigma L(2,9).$ Only $H_{1}$ has a regular set. The groups fall into the chains 
$H_{1}\subset H_{2}\subset H_{5}\subset X,$ further $H_{1}\subset H_{3}\subset H_{4}\subset X$ and $H_{3}\subset H_{5}.$

To define $H_{4}=\G(R_{4})$ select an orbit $R_{4}\subset \gO^{\{5\}}$ of length $36,$ as $X$ has $\gO^{\{5\}}$-orbits of length $180$ and $72.$ For $H_{5}=\G(R_{5})$ select an orbit $R_{5}\subset \gO^{\{5\}}$ of length $90.$ 
To define $H_{3}=\G(R_{3})$ select  $R_{3}\subset \gO^{\{5\}}$ as the union of two orbits of length $90$ and $36.$ For $H_{2}=\G(R_{2})$ select an  orbit $R_{2}\subset \gO^{\{4\}}$  of length $5.$ Finally $H_{1}=\G(R_{2}\cup R_{3}).$

\item Degree  $n=11.$ \,Here only the chain $H_{1}\subset H_{2}\subset H_{3}\subset X=M_{11}$ occurs, where $H_{1}=11,$ $H_{2}=11\e 5$ and $H_{3}=L_{2}(11).$  To define $H_{3}=\G(R_{3})$ select an orbit $R_{3}\subset \gO^{\{3\}}$ of length $55.$ The remaining two groups are contained in the relation groups $AGL(1,11)$ which has a regular set. Now apply the Basic Lemma.  

\item Degree  $n=12.$ \,We have $H_{1}=M_{11}\subset X_{1}=M_{12}$ and $H_{2}=L_{2}(11)\subset X_{2}=L_{2}(11).2.$ For $H_{1}=\G(R_{1})$  select an orbit $R_{1}\subset \gO^{\{5\}}$ of  length $132,$ and for $H_{2}=\G(R_{2})$ select an orbit $R_{2}\subset \gO^{\{6\}}$ of length $132.$ 

\item Degree  $n=13.$ \,All primitive subgroups of  $X=L_{3}(3)$ are contained in the relation group $AGL(1,13)$ which does have a regular set. Now apply the Basic Lemma.

\item Degree  $n=14.$ \,Here $H=PSL(2,13)\subset X=PGL(2,13)$ and $H=\G(R)$ where $R$ is an $H$-orbit on  $\gO^{\{3\}}.$ 

\item Degree  $n=15.$ Here $H_{1},\,H_{2},\,H_{3}\subset X=PGL(4,2)$ where $H_{1}=A_{7}$ (on the 15 cosets of $L_{3}(2)$) and $H_{2}=A_{6}\subset H_{3}=S_{6}$ (on the $2$-subsets from a $6$-set), with $H_{2}\subset H_{1}.$
For $H_{1}=\G(R_{1})$  select an orbit $R_{1}\subset \gO^{\{5\}}$ of  length $42,$ and for $H_{3}=\G(R_{3})$ select an orbit $R_{3}\subset \gO^{\{10\}}$ of length $6.$ Now $H_{2}=H_{1}\cap H_{3}=\G(R_{1}\cup R_{3}).$

\item Degree  $n=16.$ The only option is $X=AGL(4,2)=2^{4}\e L_{4}(2),$ the group of the planarity relation in ${\rm GF}_{2}^{4}.$ 
It  has $19$ primitive subgroups, four of these without regular set. We will now outline the proof for the fact that these are all groups of relations. (In brackets we state the GAP name of the group.)

There are three maximal subgroups of $X,$ these are $H_{17}=2^{4}\e S_{6}(=G_{16}),\,$ $H_{18}=2^{4}\e A_{7}(=G_{20})$ and $H_{19}=A\Gamma L(2,4)(=G_{12}).$ In each case we name an orbit of the group on subsets that is not an orbit of $X:$ 
For $H_{17}=\G(R_{17})$  select an orbit $R_{17}\subset \gO^{\{3\}}$ of  length $240,$ 
for $H_{18}=\G(R_{18})$  select an orbit $R_{18}\subset \gO^{\{5\}}$ of  length $672$ and for  $H_{19}=\G(R_{19})$  select an orbit $R_{19}\subset \gO^{\{4\}}$ of  length $20.$

We turn to groups that are maximal in one of the above. These are $H_{11}=A\Gamma L(1,16)(=G_{9}),\,$
$H_{12}=(S_{4}\times S_{4})\e 2(=G_{10}),\,$
$H_{13}=2^{4}\e S_{5}(=G_{18}),\,$
$H_{14}=2^{4}\e A_{6}(=G_{17}),\,$
$H_{15}=ASL(2,4)\e 2(=G_{13}),\,$
$H_{16}=AGL(2,4)(=G_{14}).$
The containments are $H_{12},\,H_{13},\,H_{14},\,H_{15}\,\subset H_{17},$ further  $ H_{14},\,H_{15}\,\subset H_{18}$ and $H_{11},\,H_{15},\,H_{16},\,\subset H_{19}.$
 In each case we name an  orbit of the group that is not an orbit of any group properly containing it.

\medskip
For $H_{11}=\G(R_{11}\cup R_{19} )\subset H_{19}$  select an $H_{11}$-orbit $R_{12}\subset \gO^{\{5\}}$ of  length $480.$ 
For $H_{12}=\G(R_{12}\cup R_{17} )\subset H_{17}$  select an $H_{12}$-orbit $R_{12}\subset \gO^{\{4\}}$ of  length $8.$ 
For $H_{13}=\G(R_{13}\cup R_{17})\subset H_{17}$  select an $H_{13}$-orbit $R_{13}\subset \gO^{\{4\}}$ of  length $40.$ 
We have $H_{14}=H_{17}\cap H_{18}=\G(R_{17}\cup R_{18})$  and $H_{15}=H_{18}\cap H_{19}=\G(R_{18}\cup R_{19}).$ 
For $H_{16}=\G(R_{16}\cup R_{19})\subset H_{19}$  select an $H_{16}$-orbit $R_{16}\subset \gO^{\{6\}}$ of  length $1440.$ 

The remaining groups are contained in $H_{12},\,H_{13},\,H_{15}$ or $H_{16}.$ Each of these  has a regular set of size $6,\,7$ or $8.$ These set sizes  are distinct from the set sizes of the relations defining $H_{12},\,H_{13},\,H_{15}$ and $H_{16}.$ Hence the result follows from the Basic Lemma.  

\item Degree  $n=17.$ Here we consider $H_{1}=PGL(2,2^{4})\subset H_{2}=PGL(2,2^{4}).2\subset X=P\Gamma L(2,2^{4}).$ All other primitive subgroups are contained in $AGL(1,17)$ which does have a regular set. 
Here $H_{1}$ has a regular set but $H_{2}$ does not.  To have  $H_{1}=\G(R_{1})$ we find  an $H_{1}$-orbit $R_{1}\subset \gO^{\{7\}}$  that is not an $H_{2}$-orbit.  


\item Degree  $n=21.$ We need to consider $H_{1}=PSL(3,2^{2})\subset H_{2}=P\Sigma L(3,2^{2}), \,H_{3}=PGL(3,2^{2})\subset X=P\Gamma L(3,2^{2}).$ 
To represent $H_{2}=\G(R_{2})$ we may select an $H_{2}$-orbit $R_{2}\subset \gO^{\{10\}}$  which is not an $X$-orbit and likewise,  to represent $H_{3}=\G(R_{3})$ select an $H_{3}$-orbit $R_{3}\subset \gO^{\{11\}}.$ Now $H_{1}=H_{2}\cap H_{3}=\G(R_{2}\cup R_{3}).$

 
\item Degree  $n=22.$ Here $H=M(22)  \subset X=M(22).2$  and $H$ is the automorphism group of the Steiner system on $22$ points. 

\item Degree  $n=23.$ Here $H=23\e 11 \subset X=M(23).$ Here $X$   is the automorphism group of the Steiner system on $23$ points and 
 $H$ is also a subgroup of $AGL(1,23).$  This is a relation group  and has a regular set. So we may apply the Basic Lemma. 
 
\item Degree  $n=24.$ Here $H=PSL(2,23)\subset X=M_{24}.$ In this case  $H$ is also a subgroup of $PGL(2,23)$ which is a relation group with a regular set. 

\item Degree  $n=32.$ Here we consider $H_{1}=AGL(1,2^{5})\subset H_{2}=A\Gamma L(1,2^{5})\subset X=ASL(5,2).$ To represent $H_{2}=\gO(R_{2})$ select an $H_{2}$-orbit $R_{2}\subset \gO^{\{15\}}.$  This group has a regular set of size $3$ and so $H_{1}$ is a relation group by the Basic Lemma. \dne
\end{itemize}}

\medskip

\section{ Wreath products}

Let $K$ be a permutation group on the set $\gD$ and let $L$ be a transitive permutation group on the set $\gS.$   We are interested in the imprimitive action of $K\wr L$  on $\gO:=\gD\times \gS.$  (This will be the only action considered here for wreath products.) First we 
 show that the full wreath product $(K\wr L,\,\gO) $ is a relation group if the same is true for $(K,\gD)$ and $(L,\gS)$ when $(K,\,\gD)$ is primitive.  Imprimitive groups more generally  will be considered in the next section. If $R$ is a relation then the maximum set size in $R$ is  denoted  by $m(R):=\max\{\,|x|\,:\,x\in R\}.$



\begin{thm}\label{thm:WreathProduct1} 
\,Let $\gD$ and $\gS$ be sets and put $\gO:=\gD\times \gS.$ Assume that $(K,\gD)$ and  $(L,\gS)$ are permutation groups where  $K$ acts primitively.

(i) \, If $(K,\,\gD)$ and  $(L,\,\gS)$ are relation groups then $(K\wr L,\,\gO)$ is a relation group.\\
(ii)  \,If $(K\wr L,\,\gO)$ is a relation group then $(K,\,\gD)$ is a relation group.

In the first case, if $R^{\gS}\neq \emptyset$ is a relation for $L$ then  there exists  a relation $R^{\gO}$ for  $K\wr L$  with  $m(R^{\gO})=|\gD|\,m(R^{\gS}).$
\end{thm}

\begin{cor}\label{cor:WreathProduct1} 
\,Let $(K,\gD)$ be a primitive permutation group  and let $(L,\gS)$ be  a relation group. Then   $K\wr L$ is a relation group on $\gD\times \gS$ if and only if $(K,\gD)$ is not $ \Al(\gD)$ and  not one of the groups in the list ${\cal L}$ of Section 4.
\end{cor}

{\sc Comment:\,} It is essential here that we consider the full wreath product. As an illustration take $(K,\gD)=(C_{2},\{1,2\})=(L,\gS)$ when $C_{2}\wr C_{2}=D_{8}$ is a relation group on $\{1,2,3,4\}.$ The Klein group $V\subset C_{2}\wr C_{2}$ is orbit closed, as we have seen earlier, but it is not a relation group. Note that $D_{8}$ has no regular set. This is not an isolated example:

\begin{cor}\label{cor:WreathProduct11} 
\, There are infinitely many wreath products $(K\wr L,\,\gO)$ with both $(K,\,\gD)$ and  $(L,\,\gS)$ primitive such that $(K\wr L,\,\gO)$ is orbit closed but not a relation group. 
\end{cor}



\medskip
{\it Proof of Theorem~\ref{thm:WreathProduct1}:} \, For $(i)$\, let $\gS=\{1..s\}$ and denote $\gD_{i}:=\gD\times \{i\}$ where we identify  $\gD=\gD_{1}.$ Suppose that there is a relation  $R^{\gD}$ on $\gD$ so that $K=\G(R^{\gD}).$ Since $K$ is transitive we may assume that all sets $y$ in $R^{\gD}$ have size $2\leq |y|\leq d-2.$  Correspondingly let $R^{\gD_{i}}\subseteq P(\gD_{i})$ be a copy of $R^{\gD}$ on $\gD_{i}.$ Let \begin{equation}\label{rela1}R'=\big\{\,\gO\setminus y\,\,: \,\,\mbox{$y$ belongs to some $R^{\gD_{i}}$}\big\}\,.\end{equation}

Let $H=K\wr \Sy(\gS).$ We claim that  $H=\G(R').$ Clearly, $H\subseteq \G(R').$ On $\gO$ we define an adjacency relation   $\alpha\sim\beta$ if and only if there is some $x\in R'$ with  $|\{\alpha,\,\beta\}\cap x|=0.$ (In other words, $\alpha\sim\beta$ if and only if  $\{\ga,\gb\}\subseteq y$ for some $y\in R^{\gD_{i}}.)$
Evidently there is no edge from $\gD_{i}$ to $\gD_{j}$ if $i\neq j,$ and as $\G(R^{\gD_{i}})$ is primitive $\gD_{i}$ is a connected component of the graph, for $i=1,..,s.$ Therefore each element $g$ in $\G(R')$ permutes the $\gD_{i}$ as sets, and if $g$ stabilizes $\gD_{j}$ then it preserves $R^{\gD_{j}}.$ Hence $g\in H$ and therefore  $\G(R')=H.$ Now 
\begin{equation}\label{rela2}R=R^{\gD_{1}}\cup ...\cup R^{\gD_{s}}\,\end{equation}
is the collection of all complements of sets in $R'$ and so also $H=\G(R).$ All sets in $R$ have size $\leq d-2.$

Next let $R^{\gS}$ be a relation of $\gS$ with $L=\G(R^{\gS})$ and suppose that the maximum set size in $R^{\gS}$ is $m_{\gS}.$ Define \begin{equation}\label{rela3}R''=\big\{\,\gD_{i_{1}}\cup...\cup \gD_{i_{t}}\,\,:\,\,\{i_{1}...i_{t}\}\in R^{\gS}\,\big\}\,.\end{equation}
The sets in $R''$ have size $\ell d$ where $1\leq \ell\leq m_{\gS}$ and they are larger  than any set in $R.$ It follows that   $$\G(R\cup R'')=\G(R)\cap \G(R'').$$

In particular, $K\wr L\subseteq \G(R\cup R'').$ Conversely, if $h\in \G(R\cup R'')$ then $h\in \G(R)=H,$ and so $h$ permutes the sets $\gD_{i}.$ Therefore $h$ induces a permutation $h_{\gS}$ of $\gS$ preserving $R^{\gS}.$ Hence $h_{\gS}\in L=\G(R^{\gS})$ and so $h$ belongs to $K\wr L.$ Therefore $K\wr L=\G(R\cup R'')$ is a relation group on $\gO.$ 

\medskip
Conversely, for $(ii)$\, suppose that there is some relation $R$ on $\gO$ so that 
$K\wr L=\G(R).$ For each $i=1...s$ consider the relation 
$R_i:=\{ y\cap \gD_i\,\,:\,\, y\in R \}$ and the subgroup of $K\wr L$ given by  
$K_i=\{(1,\dots ,1,k,1,\dots,1)\,\,:\,\, k\in K \}$ with $k$ in position $i.$
Then $K_{i}$ fixes  $\gD_i$ and preserves $R_{i}$ so that $ K_{i}\subseteq \G(R_{i}).$ 
As $L$ is transitive the permutation groups $(K_{i},\,\gD_{i})$ are permutationally isomorphic to each other for $i=1...s.$ In particular,   $K=K_{1}\subseteq \G(R_{1})$ as a group on $\gD=\gD_{1}.$

Assume that $K$ is not a relation group. From Corollary~\ref{cor:PrimRela}\, it follows  that $K$ is not orbit closed. Let therefore  $K^{*}\neq K$ be its orbit closure, so that  $K\subset K^{*}\subseteq G(R_{1})\neq K$ by Proposition~\ref{2.5}.  We show that this implies $ K^{*}\wr L\subseteq \G(R),$ contradicting that $K\wr L= \G(R).$ 

A subset $x$ of $\gO$ is of the form $x=\bigcup_{i=1...s}\, (x_{i},i)$  with $x_{i}\subseteq \gD$  and for $h=(k_{1},..,k_{s};b)\in K\wr L$ we have $xh=(x_{1}k_{1},1b)\cup...\cup (x_{s}k_{s},sb).$
To show that  $ K^{*}\wr L\subseteq \G(R)$ let $h^{*}=(k^{*}_{1},..,k^{*}_{s};b)\in K^{*}\wr L$ and let $x$ be a set  in $R.$   Then by the orbit equivalence of $K$ and $K^{*}$ there are elements $k_{1},..,k_{s}$ in $K$ so that $x_{i}k^{*}_{i}=x_{i}k_{i}$ for $i=1...s.$ Let $h=(k_{1},..,k_{s};b)\in K\wr L$ for this choice of the $k_{i}.$ Then 
$xh^{*}=(x_{1}k^{*}_{1},1b)\cup...\cup (x_{s}k^{*}_{s},sb)=(x_{1}k_{1},1b)\cup...\cup (x_{s}k_{s},sb)=xh$ belongs to $R.$ Hence $K$ is a relation group.\dne

{\it Proof of Corollaries~\ref{cor:WreathProduct1} and ~\ref{cor:WreathProduct11}:} \, The first corollary is immediate. In Seress and Yang~\cite{Seress2}\, the orbit closed groups of type $(K\wr L,\,\gO)$ with both $(K,\,\gD)$ and  $(L,\,\gS)$ primitive are classified, when $(K,\,\gD)$ has a regular set and  $L$ does not contain $\Al(\gS).$ Among these groups one may take for instance  $(K,\,\gD)=(C_{3},\,3)$ and any $L$ with $|\gS|>9.$  As $(C_{3},\,3)$ is not a relation group  Theorem~\ref{thm:WreathProduct1}\, implies that  $(K\wr L,\,\gO)$ is not a relation group. \dne

\bigskip 
Next we are interested in conditions which guarantee that  $K\wr L$ has a regular set provided   $(K,\gD)$ and $(L,\gS)$ have suitable regular sets.  

\begin{prop}\label{prop:RelGr2}
 \,Let $\gD$ be a set of size $d$ and  let $\gS$ be a set of size $s.$  Assume that $(K,\gD)$ and $(L,\gS)$ are  transitive permutation groups with regular sets of  size $r_{\gD}\neq\frac{d}{2}$ and  $r_{\gS}$, respectively. 
 Then $K\wr L$ has a regular set of size $r_{\gO}=r_{\gS}d+(s-2r_{\gS})(d-r_{\gD})$ as an imprimitive  group on $\gO=\gD\times \gS.$ \end{prop}

\pf  Let $x\subset \gD$ be a regular set for $(K,\gD)$ of size $r_{\gD}\neq \frac{d}{2}$  and let $\{1,\,..,\,k\}\subset \gS$  be a regular set  for $(L,\gS)$ of size with $k=r_{\gS}.$ Define $x'=\gD\setminus x$ and consider the following subset $w$ of $\gO=(\gD\times \{1\})\,\cup ... \,\cup \,(\gD\times \{s\})$ defined by \begin{equation}\label{rela4}w=(x\times \{1\})\cup...\cup  (x\times \{k\})\,\,\cup \,\,(x'\times \{k+1\})\cup...\cup (x'\times \{s\})\,\,.\end{equation} Then we have $|w|=r_{\gS}d+(s-2r_{\gS})(d-r_{\gD}).$

We claim that $w$ is a regular sets for $K\wr L.$ Let $h\in K\wr L$ with $w^{h}=w.$ Since $h$ preserves the blocks and since $|x|\neq |x'|$ it follows that $h$ permutes the block sets $\big\{(x\times \{1\})\cup...\cup  (x\times \{k\})\big\}$  and $\big\{(x'\times \{k+1\})\cup...\cup (x'\times \{s\})\big\}$ separately. Since $\{1..k\}$ and $\{k+1,..,s\}$ are regular sets for $L$ we have that $(x\times\{i\})h=(x\times\{i\})$ for $1\leq i\leq k$ and $(x'\times\{j\})h=(x'\times\{j\})$ for $k+1\leq j\leq s.$ But as $(x\times\{i\})$ and $(x'\times\{j\})$ are regular sets for the $H$ action on $(\gD\times\{i\})$ and $(\gD\times\{j\}),$ respectively, we have that $h$ is the identity on $\gD\times \{1..s\}.$ \dne

By a similar argument we can show that $K\wr L$ has a regular set provided $K$ has sufficiently many sets in comparison to the degree of $L.$

\begin{prop}\label{prop:RelGr3} \,Let $(K,\gD)$ and  $(L,\gS)$ be permutation groups where $L$ is transitive.  Assume that $K$ has at least $|\gS|$ regular sets on $\gD$ of pairwise distinct  size.  Then $K\wr L$ has a regular set as an imprimitive group on $\gO=\gD\times \gS.$ \end{prop}

\pf  Let $s:=|\gS|$ and let $x_{1},..,\,x_{s}\subset \gD$ be  regular sets for $K$ of pairwise distinct sizes.  Then \begin{equation}\label{rela15}w=(x_{1}\times \{1\})\cup...\cup (x_{s}\times \{s\})\,\,\end{equation} 
is a regular set for $K\wr L$ by the same reasoning as in the proof of Proposition~\ref{prop:RelGr2}.  \dne

{\sc Comment:}\, Trivially, if a group of degree $n$ admits regular sets then it has regular sets of size $\leq \frac n2,$ just replace a set by its complement.  In particular, in  Proposition~\ref{prop:RelGr2} one may select regular sets to have size  $r_{\gD}<\frac d2$ and $r_{\gO}< \frac{ds}2$ provided $r_{\gS}\neq \frac s2.$

\bigskip

\section{ Imprimitive Groups}

Let $H$ be   an  imprimitive permutation group on $\gO$ and let  
$\gO_{1},.., \gO_{s}$ with $\gO=\gO_{1}\cup ...\cup \gO_{s}$ be a system of blocks of imprimitivity. (As usual, an imprimitive group is transitive.) Then $H$ induces permutation groups $K:=H^{\gD}_{\gD}$ on $\gD:=\gO_{1}$ and $L:=H^{\gS}$ on  $\{\gO_{1},.., \gO_{s}\}$ where  $\gS:=\{1,..,s\}.$ In particular, $H\subseteq K\wr L$ acts in the imprimitive action on $\gO=\gD\times \gS.$ Evidently $K$ is primitive on $\gD$ if and only if  $\gD$ is a minimal (non-trivial) block of imprimitivity of $H.$ If $R$ is a relation recall that ${\rm ar}(R)=\{ |x|\,:\,x\in R
\}$ is the arity of $R$ and the maximum set size in $R$ is  $m(R):=\max\{\,|x|\,:\,x\in R\}.$


\begin{thm} \label{thm:many} 
\, Let $H$ be an imprimitive permutation group on $\gO$ and let $\gD\subset \gO$ be a minimal block of imprimitivity. Suppose that $H^{\gD}_{\gD}$  is a relation group on $\gD$ which has at least $\frac{|\gO|}{|\gD|}$ regular sets of pairwise distinct size. Then every subgroup of $H$ is a relation group.  
\end{thm}

\pf We put  $K:=H^{\gD}_{\gD},$ $d:=|\gD|$ and $n:=|\gO|.$ 
Let $R^{\gD}$ be a relation on $\gD$ with $K=\G(R^{\gD}).$ Let $\gS$ be a set of size $s:= \frac{n}{d},$ hence  $\gO=\gD\times \gS,$ and let $R^{\gD_{i}}$ be a copy of $R^{\gD}$ on $\gD\times \{i\}.$
Then $H\subseteq \bar H :=K\wr \Sy(\gS)$ and $\bar H=\G(R')$ where $R'$ is the  relation defined in \,(\ref{rela1})\, in the proof of   Theorem~\ref{thm:WreathProduct1}. All sets in $R'$ have size $\geq (s-1)d+2>\frac {n}{2}$ since $s\geq 2.$ It follows from  Proposition~\ref{prop:RelGr3}\, that $\bar H$ has a regular set of size $\leq  \frac {n}{2}.$ Now apply the first part of the Basic Lemma. \dne

Asymptotically most sets in a primitive permutation not containing the alternating group of the same degree are regular. More precisely we have the following:

\begin{cor} \label{cor:many} 
\, There is a function $d(s)$ so that the following is true. Let $H$ be an imprimitive permutation group of degree $n$ with a block $\gD$ of imprimitivity so that  \\
(i) \,\, $H^{\gD}_{\gD}$ does not contain  $\Al(\gD),$ \\
(ii) \, $|\gD|\geq d(s)$ for  $s=\frac{n}{|\gD|}.$\\
Then every  subgroup of $H$ is a relation group. 
\end{cor}
 
\pf We sketch the argument.  Let  $s$ be fixed.  Then there is a  bound $d(s)$  so that any primitive group $(K,\,\gD)$ not containing $\Al(\gD)$ with $|\gD|\geq d(s)$ has at least $s$ regular sets on $\gD$ of pairwise different size.  This is shown in Theorem~\ref{t=2}\, in Section 7\, for $s=2$ and this argument extends readily to arbitrary $s.$ Now apply Theorem~\ref{thm:many}.\dne





\begin{thm}\label{RelGr3} 
Let $H$ be an imprimitive permutation group on $\gO$ with a minimal block $\gD\subset \gO$ of imprimitivity. Put  $K:=H^{\gD}_{\gD},$ $d:=|\gD|$ and $n:=|\gO|.$ Let $(L,\gS)$ with $s:=|\gS|=\frac nd$ be the group induced by $H$ on the blocks of imprimitivity. Suppose that \vspace{-5mm} \begin{enumerate}[(i)]
\item $K$ is a primitive relation group and has a regular set of size $r_{\gD}\neq\frac{d}{2}$ on $\gD,$ and\vspace{-2.5mm}
\item $L$ is a relation group  and  $L$ has a regular set on $\gS.$ 
 \vspace{-5mm}
\end{enumerate}

Then every subgroup of $H$ is a relation group on $\gO.$
\end{thm}

\pf  We follow the proof of Theorem~\ref{thm:WreathProduct1}\, using  the same terms. Let $R$ and $R''$ be as in \,(\ref{rela2})\, and \,(\ref{rela3}), and let $w$ be the regular set for $K\wr L$ as in \,(\ref{rela4}). Let $G\subseteq K\wr L.$ We put $R^{G}= w^{G},$ the orbit of $w$ under $G,$  and consider 
\begin{equation}\label{rela5} G':=\G(R\cup R''\cup R^{G})\,.\end{equation}

Clearly $G\subseteq G'$ and it remains to show the converse. First note that a set in $R$ has size $\leq d-2$ while sets in $R''$ have size $\ell d$ with $1\leq \ell\leq m_{\gS}.$ The set $w$ has size $> d$. Therefore $R$ is preserved by $G'$ as a set and  so 
\begin{equation}\label{rela6} G'=\G(R)\cap \G(R''\cup R^{G})\subseteq K\wr\Sy(\gS)\,,\end{equation} by  the  proof above of Theorem~\ref{thm:WreathProduct1}. Let $h$ be an element of $G'.$ If $h$ leaves $R^{G}$ invariant then also $R''$ is invariant under $h$ so that $h\in \G(R\cup R'')=K\wr L.$ This means that there is some $g\in G$ with  $w^{h}=w^{g}$ so that $hg^{-1}\in K\wr L$ stabilizes $w.$ Since $w$ is a regular set for $K\wr L$ we have $h=g\in G.$

Hence we have to consider the situation  where $h$ does not leave $R^{G}$ invariant. We claim that this is impossible. Let $u\in R^{G}$ with $u^{h}\in R''.$ Since $G\subseteq K\wr L$ it follows that $u$ decomposes over the $\gD_{i}$ in the same way \,(\ref{rela4}) as $w$ so that  $$u=u_{1}\cup...\cup u_{k}\cup u'_{k+1}\cup...\cup u'_{s}$$
where $2\leq |u_{1}|=...=|u_{k}|=r_{\gD}\leq d-2$ and  $d-r_{\gD}=|u_{k+1}|=...=|u_{s}|.$ Since $h\in K\wr \Sy(\gS)$ preserves the block structure also $u^{h}$ decomposes in this way. But as $u^{h}\in R''$ this is in contradiction to \,(\ref{rela3}). \dne


\bigskip

\section{ Some Classes of Imprimitive Relation Groups}

Let  $H$ be   an  imprimitive permutation group on $\gO$ and let  $\gD^{(0)}$ be some minimal (non-trivial) block of imprimitivity of $H.$ Let $K^{(0)}$ be the group induced by $H$ on $\gD^{(0)}$ and let $L$ be the group induced by $H$ on the $H$-images of $\gD^{(0)}.$ Then $H\subseteq K^{(0)}\wr L.$ 

By   repeating the same process for $L$ we obtain a chain $(K^{(0)},\gD^{(0)}),..,(K^{(t)},\gD^{(t)})$ of primitive groups so that $H$ is contained in  the iterated wreath product, 
$$H\subseteq K^{(0)}\wr(K^{(1)}\wr(...\wr(K^{(t-1)}\wr K^{(t)})...))\,.$$ 
We call $(K^{(0)},\gD^{(0)}),..,(K^{(t)},\gD^{(t)})$ with $t\geq 1$ a {\it imprimitivity chain of length $t$} for $H.$ As there is no Jordan-H\"older Theorem for primitivity a group may  have several  incomparable imprimitivity chains, including chains of different lengths. 

We say that $H$ is  {\it ${\cal A}\!$-imprimitive } if there is a collection ${\cal A}$  of primitive permutation groups  and some  imprimitivity chain $(K^{(0)},\gD^{(0)}),..,(K^{(t)},\gD^{(t)})$ for $H$ so that $(K^{(i)},\gD^{(i)})$ belongs to ${\cal A}$ for all $i=0...t.$ Similarly, $H$ is  {\it ${\cal A}'\!$-imprimitive } if there is an imprimitivity chain $(K^{(0)},\gD^{(0)}),..,(K^{(t)},\gD^{(t)})$ for $H$ so that $(K^{(i)},\gD^{(i)})$ does not belong to ${\cal A}$ for all $i=0...t.$
(In both cases the containment is  up to permutation isomorphisms.) Note that a group may happen to be  ${\cal A}\!$-imprimitive  and ${\cal A}'\!$-imprimitive  at the same time.

In view of Theorem~\ref{RelGr3} we are interested in those   primitive groups which have regular sets, all  of the same cardinality, 
$$\lur:=\big\{\,\, \mbox{$(K,\gD)$ is primitive with regular sets, all  of size $\frac 12 |\gD|$}\,\,\big\}.\,$$

Now let 
\begin{eqnarray}{\cal A}&:=&
\big\{\,\Sy(\gD)\, :\,\,|\gD|\geq 2\,\big\}\,\cup \, \big\{\, \Al(\gD)\,:\,\,|\gD|\geq 3\,\big\}\nonumber\\
&\cup&\lnr \,\,\,\cup\,\,\, \lur\nonumber\\
&\cup&\big\{  (5,C_{5}),\,  (8,2^{3}\e 7), \,
(9, 3^{2}\e 8) \big\}\,.
\end{eqnarray}

It is easy to check that the list $\lst$ and $\loe$ from Section 4\, are contained in ${\cal A}.$ Now  we have the following

\begin{thm}\label{inductionGluck} 
 \, Let $H$ be an ${\cal A}'\!$-imprimitive group on $\gO.$  Then all subgroups of $H$ are relation groups on $\gO.$ 
\end{thm}

\pf Let $H\subseteq K^{(0)}\wr(K^{(1)}\wr(...\wr(K^{(t-1)}\wr K^{(t)})...))$
be an imprimitivity chain for $H$ where $(K^{(i),\,\gD^{(i)}})$ is not in  ${\cal A},$ for all $i=0...t.$ As the list ${\cal L}$ from Section 4\, is contained in ${\cal A}$  it follows from Theorem~\ref{thm:PrimRela}\,  that $(K^{(i)},\gD^{(i)})$ is a relation groups for all $i=0...t.$ In addition, all $(K^{(i)},\gD^{(i)})$ have regular sets of cardinality $\neq \frac 12 |\gD^{(i)}|.$ 
It follows by induction from  Proposition~\ref{prop:RelGr2}\, that $(K^{(1)}\wr(...\wr(K^{(t-1)}\wr K^{(t)})...))$ is a relation group with a regular set, see the comment following Proposition~\ref{prop:RelGr3}. 
Now use Theorem~\ref{RelGr3}. \dne

Before we turn to the groups in $\lur$ specifically we are able to state several corollaries which do not  require the explicit knowledge of this class of groups. 

\bigskip
{\sc Odd Degree and Odd Order:}\, A primitive group of odd degree with at least one  regular set automatically has two regular sets of different sizes.  Therefore 

\begin{cor}  Let $H$ be an ${\cal O}'\!\!$-imprimitive group of odd degree where ${\cal O}=\{\Sy(\gD),\, \Al(\gD)\,:\,\,|\gD|\geq 3\}\,\cup\,\lnr \,\cup\,\{ (5,C_{5}),\
(9, 3^{2}\e 8) \}.$ Then all subgroups of $H$ are relation groups. 
\end{cor}

In particular, 
\begin{cor}  Let $H$ be an ${\cal O}'\!\!$-imprimitive group of odd order where ${\cal O}=\{(3,C_{3}),\, (5,C_{5})\}.$ Then all subgroups of $H$ are relation groups. 
\end{cor}

\medskip
{\sc Solvable Groups:}\, The regular sets of primitive solvable groups have been studied in Gluck~\cite{gluck}.  From his result it follows that the solvable groups belonging to  ${\cal A}$ form the list  
\begin{eqnarray}{\cal S}&:=&
\big\{\,\Sy(\gD),\,\Al(\gD)\,\, :\,\,|\gD|\leq 4\,\big\}
\nonumber\\
&\cup&    
\big\{(5,D_{10}),\,(7,F_{42}),\,(8,2^{3}\e 7.3),\,
(9,3^{2}\e (2.L_{2}(3)),\,(9,3^{2}\e 2.L_{2}(3).2),\,
(9,3^{2}\e D_{8}),\,(9,3^{2}\e 8.2)\big\}\nonumber\\
&\cup&\big\{  (5,C_{5}),\,  (8,2^{3}\e 7), \,
(9, 3^{2}\e 8) \big\}\,.
\end{eqnarray}
Therefore the theorem provides the following

\begin{cor}  Let $H$ be a solvable  ${\cal S}'\!\!$-imprimitive group. Then all subgroups of $H$ are relation groups. 
\end{cor}

Using much the same idea as in \cite{CNS} we show that $\lur$ is finite, and this will conclude the  discussion of this particular class of groups. In a forthcoming paper we  will determine $\lur$ explicitly. 

\begin{thm} \label{t=2}
Let $G$ be a primitive permutation group on $\gO$ and assume that $G$ does  not contain $\Al(\gO).$ Then with only finitely many exceptions $G$ has regular sets of cardinality $k_{1}\neq k_{2}$. 
\end{thm}
 
\pf  This follows from Theorem~\ref{Seress1}\, if $n=|\gO|$ is odd. So let $n=2k$ and denote the minimal degree of $G$ by $m,$ so that $m$ is  the least  number of elements in $\gO$ moved by a non-identity element of $G.$ Hence a non-identity element fixes at most $2^{n-\frac m2}$ subsets of $\gO.$ Therefore the cardinality of the set $$A=\{(\gD,g)\,\,:\,\,|\gD|\neq k,\,\,\gD^{g}=\gD,\,\,\,g\neq 1\}$$ can be bounded  by  
$$ 2^{n}-{n\choose k}\,\,\leq\,\, |A|\,\,\leq \,\, 2^{n-\frac m2}\,(|G|-1) .$$ 
As the number of odd-sized (or even-sized) subsets of an $n$-element set is $2^{n-1}$ we have 
$$ 2^{n-1}\,\,\leq\,\,2^{n}-{n\choose k}\,\,\leq 2^{n-\frac m2}\,|G|.$$
In particular, $2^{\frac m2-1}\,\,\leq\,\, |G|,$ and from the well-known bound
$2^{\frac n{m}}\leq |G|,$ see for instance the Proposition in~\cite{CNS}, we have 
$$\frac n2 \leq \log|G| \cdot (\log|G| +1).$$ In particular $$n^{1+(\log_2n)}\,\,<\,\,|G|$$ for all sufficiently large $n.$ Using Mar\'oti's estimate in \cite{Maroti}\, for the order of a primitive groups it follows that $G$ is a subgroup of $\Sy_{\!s}\wr \Sy_{\!t}$ acting in the product action of degree ${s \choose a}^{t}$ where $\Sy_{\!s}$ acts naturally on the $a$-element subsets of $\{1...s\},$ for some $a$ and $t\geq 2.$ In Lemma~3.1 of \cite{Seress2}\, it is shown that in this action $\Sy_{\!s}\wr \Sy_{\!t}$ has regular sets of at least four distinct sizes. \dne

\bigskip

\section{Closing Remarks and  Open Problems}

We conclude with a few remarks and suggestions for further work on the two  topics of this paper. 

{\sc Relation Groups:}\, From the results in Sections 5--7 one should take the view that a transitive group in general tends to be a relation group, and that it may not be so easy to classify the imprimitive groups that are not relation groups. The only firm result to show that certain kinds of imprimitive groups are not relation groups is Theorem~\ref{thm:WreathProduct1}. A measure for a group $(G,\,\gO)$  to be a relation group is 
$$r(G):=\min\{\,|\G(R):G|\,\,\,\,\mbox{where $R$ is a relation on $\gO$ with $G\subseteq \G(R)$}\}.$$

Thus  $r(G)=1$ if and only if $G$ is a relation group.

\begin{itemize}\vspace{-3mm}
\item Is  $r(G)$ bounded in terms of the length of an imprimitivity chain for $G?$ Are there absolute bounds for good  classes of imprimitive group?
\item In all examples with  $r(G)>1$ we are aware of it appears that $r(G)$ is even. Is this true for all imprimitive groups? 
\end{itemize}

{\sc Orbit Closure:}\,  In the literature orbit closure properties have been considered almost exclusively for  groups with imprimitivity chains of length $\leq 2.$ (Groups with length $2$ are the two-step imprimitive groups in \cite{Seress2}.)  The results here on relation groups imply closure properties for arbitrary subgroups in the presence of suitable regular sets.   A measure for a groups $(G,\,\gO)$ to be orbit  closed is
$$c(G):=|G^{*}:G|$$
where $G^{*}$ is the orbit closure of $G.$ By Proposition~\ref{2.5}\, we have $G\subseteq G^{*}\subseteq \G(R)$ for any relation $R$ with $G\subseteq \G(R)$ and therefore $c(G)$ divides $r(G).$

\begin{itemize}\vspace{-3mm}
\item The set-transitive group $G=L_{2}(8)$ of degree $9$ in Lemma~\ref{2.1}\, has $c(G)=6!.$ This is the largest value for $c(G)$ we are aware of. Is there an absolute bound for $c(G)?$
\item Considering Corollary~\ref{cor:WreathProduct11}, what are the groups $G$ with $1=c(G)<r(G)?$   
\end{itemize}

Finally note that if $G\subseteq \G(R)$ for a relation $R$ then $G\subseteq G^{*}\subseteq \G(R)$ 
and hence
$r(G)=c(G)\cdot r(G^{*}).$ This suggests that questions about relation groups can  be reduced to orbit closed groups.

\end{document}